\documentclass[a4paper, 12pt]{article}
\usepackage[top = 2.5cm, bottom = 2.5cm, left = 2.5cm, right = 2.5cm]{geometry}
\usepackage{cite, amsmath, amssymb}
\usepackage[margin = 1cm, font = small, format = hang, labelsep = period, labelfont = bf]{caption}
\usepackage[latin1]{inputenc}
\usepackage{amsfonts}
\usepackage{amsthm}
\usepackage{makeidx}
\usepackage{graphicx}
\usepackage{mathrsfs, xcolor}
\usepackage{enumerate}
\usepackage{blkarray, bigstrut}
\usepackage{bm}
\usepackage{setspace}
\usepackage{float}
\usepackage[ruled, vlined]{algorithm2e}
\usepackage{lineno}
\usepackage{etoolbox}
\usepackage{authblk}
\usepackage{hyperref} 

\numberwithin{table}{section}
\numberwithin{equation}{section}
\theoremstyle{plain}
\newtheorem{theorem}{Theorem}[section]
\newtheorem{proposition}[theorem]{Proposition}
\newtheorem{definition}[theorem]{Definition}

\newtheorem{example}[theorem]{Example}

\newtheorem{corollary}[theorem]{Corollary}
\newtheorem{remark}[theorem]{Remark}

\renewenvironment{equation*}
	{\begin{linenomath*}\begin{equation}\nonumber}
	{\end{equation}\end{linenomath*}}
\newenvironment{eqnarray'}
	{\begin{linenomath*}\begin{eqnarray*}}
	{\end{eqnarray*}\end{linenomath*}}

\title{\vspace{3.5cm}\textbf{Concentration Robustness in \\ LP Kinetic Systems}}
\date{}

\author[1,2,3,*]{\textbf{Angelyn R. Lao}}
\author[2]{\textbf{Patrick Vincent N. Lubenia}}
\author[4]{\\ \textbf{Daryl M. Magpantay}}
\author[1,2,5,6]{\textbf{Eduardo R. Mendoza}}

\affil[1]{\small\textit{Department of Mathematics and Statistics, De La Salle University, 2401 Taft Avenue, Manila, 0922, Metro Manila, Philippines}}
\affil[2]{\small\textit{Systems and Computational Biology Research Unit, Center for Natural Sciences and Environmental Research, 2401 Taft Avenue, Manila, 0922, Metro Manila, Philippines}}
\affil[3]{\small\textit{Center for Complexity and Emerging Technologies, 2401 Taft Avenue, Manila, 0922, Metro Manila, Philippines}}
\affil[4]{\small\textit{College of Arts and Sciences, Batangas State University, Batangas, 4200, Philippines}}
\affil[5]{\small\textit{Max Planck Institute of Biochemistry, Am Klopferspitz 18, 82152, Martinsried near Munich, Germany}}
\affil[6]{\small\textit{Faculty of Physics, Ludwig Maximilian University, Geschwister-Scholl-Platz 1, 80539, Munich, Germany}}
\affil[*]{Corresponding author: \texttt{angelyn.lao@dlsu.edu.ph}}

\begin{document}
\maketitle
\thispagestyle{empty}


\begin{abstract}
    For a reaction network $\mathscr{N}$ with species set $\mathscr{S}$, a log-parametrized (LP) set is a non-empty set of the form $E(P, x^*) = \{x \in \mathbb{R}^\mathscr{S}_> \mid \log x - \log x^* \in P^\perp\}$ where $P$ (called the LP set's flux subspace) is a subspace of $\mathbb{R}^\mathscr{S}$, $x^*$ (called the LP set's reference point) is a given element of $\mathbb{R}^\mathscr{S}_>$, and $P^\perp$ (called the LP set's parameter subspace) is the orthogonal complement of $P$. A network $\mathscr{N}$ with kinetics $K$ is a positive equilibria LP (PLP) system if its set of positive equilibria is an LP set, i.e., $E_+ (\mathscr{N}, K) = E(P_E, x^*)$ where $P_E$ is the flux subspace and $x^*$ is a given positive equilibrium. Analogously, it is a complex balanced equilibria LP (CLP) system if its set of complex balanced equilibria is an LP set, i.e., $Z_+ (\mathscr{N}, K) = E(P_Z, x^*)$ where $P_Z$ is the flux subspace and $x^*$ is a given complex balanced equilibrium. An LP kinetic system is a PLP or CLP system. This paper studies concentration robustness of a species on subsets of equilibria, i.e., the invariance of the species concentration at all equilibria in the subset. We present the ``species hyperplane criterion'', a necessary and sufficient condition for absolute concentration robustness (ACR), i.e., invariance at all positive equilibria, for a species of a PLP system. An analogous criterion holds for balanced concentration robustness (BCR), i.e., invariance at all complex balanced equilibria, for species of a CLP system. These criteria also lead to interesting necessary properties of LP systems with concentration robustness. Furthermore, we show that PLP and CLP power law systems with Shinar-Feinberg reaction pairs in species $X$, i.e., their rows in the kinetic order matrix differ only in $X$, in a linkage class have ACR and BCR in $X$, respectively. This leads to a broadening of the ``low deficiency building blocks'' framework introduced by Fortun and Mendoza (2020) to include LP systems of Shinar-Feinberg type with arbitrary deficiency. Finally, we apply our results to species concentration robustness in LP systems with poly-PL kinetics, i.e., sums of power law kinetics, including a refinement of a result on evolutionary games with poly-PL payoff functions and replicator dynamics by Talabis et al (2020).
\end{abstract}


\baselineskip = 0.30in


\section{Introduction}
\label{sec:intro}

A log-parametrized (LP) set is a non-empty subset of $\mathbb{R}^\mathscr{S}_>$ (the set of positive real-valued functions with domain $\mathscr{S}$) of the form $E(P, x^*) := \{x \in \mathbb{R}^\mathscr{S}_> \mid \log x - \log x^* \in P^\perp\}$ where $P$ (called the LP set's flux subspace) is a subspace of $\mathbb{R}^\mathscr{S}$, $x^*$ (called the LP set's reference point) is a given element of $\mathbb{R}^\mathscr{S}_>$, and $P^\perp$ (called the LP set's parameter subspace) is the orthogonal complement of $P$. Feinberg, Horn, and Jackson discovered interesting relationships between positive equilibria and LP sets of the stoichiometric subspace of a mass action system (i.e., a system with mass action kinetics) in the early 1970's, which are recorded in Feinberg's 1979 Wisconsin Lecture Notes \cite{FEIN1979}. In particular, for a reaction network $\mathscr{N}$ with kinetics $K$ and which is absolutely complex balanced, i.e., those where the non-empty sets of positive equilibria $E_+ (\mathscr{N}, K)$ and complex balanced equilibria $Z_+ (\mathscr{N}, K)$ coincide, one has $E_+ (\mathscr{N}, K) = Z_+ (\mathscr{N}, K) = E(S, x^*)$ where $S$ is the stoichiometric subspace of $\mathscr{N}$ and $x^*$ is a given equilibrium. Abstracting from this pioneering work, we call a kinetic system ``of positive equilibria LP (PLP) type'' \linebreak (``of complex balanced equilibria LP (CLP) type'') if its non-empty set of positive equilibria (complex balanced equilibria) is an LP set. A system of LP type (or simply an \linebreak ``LP system'') is a system of PLP or CLP type.

Various LP systems beyond mass action systems have been studied. In 2014, S. M\"{u}ller and G. Regensburger \cite{MURE2014} showed that any complex balanced generalized mass action system (GMAS) is a CLP system whose flux subspace $\tilde{S}$ is its kinetic order subspace. GMAS include all power law systems (i.e., systems with power law kinetics) where branching reactions of any reactant complex have identical rows in the kinetic order matrix (called reactant-determined kinetic (PL-RDK) systems). In 2018, Talabis et al \cite{TAM2018} demonstrated that all systems satisfying the Deficiency One Theorem for power law systems with \linebreak $\hat{T}$-rank maximal kinetics (PL-TIK) are of PLP type. In 2019, part of the Deficiency Zero Theorem of Fortun et al \cite{FLRM2019} established that the class of non-PL-RDK systems (called PL-NDK systems) are of PLP type. A broad generalization of this result was derived by B. Hernandez and E. Mendoza \cite{HEME2021}. Subsets of poly-PL systems (i.e., sums of power law systems) and Hill-type systems have also been shown to be LP systems \cite{FTJM2020, HEME2020}.

In this paper, we study properties of concentration robustness of species in LP systems. G. Shinar and M. Feinberg \cite{SHFE2010} introduced the concept of absolute concentration robustness (ACR) in a mass action system: a species has ACR if its value at all positive equilibria of the system is the same. They presented a sufficient condition for ACR in deficiency one networks which was remarkably abstracted from subsystems in the bacterium \textit{Escherichia coli}. They called the condition ``structural'' as it was based on the occurrence of a pair of complexes in the network. In 2018, Fortun et al \cite{FMRL2018} extended their result to deficiency one PL-RDK systems by reinterpreting the properties of the pair of complexes in terms of the kinetic orders of the corresponding reaction pairs (which are called Shinar-Feinberg pairs or SF-pairs), indicating the primarily kinetic nature of the condition. In 2020, N. Fortun and E. Mendoza \cite{FOME2021} extended the SF-pair condition to deficiency zero systems using dynamic equivalence, further highlighting the property's primarily kinetic character. They also introduced the weaker concept of balanced concentration robustness (BCR) which required that species have the same value only on a subset of complex balanced equilibria. These low deficiency networks, i.e., with deficiency 0 or 1, became ``building blocks'' of concentration robustness in larger and higher deficiency power law systems in independent (for ACR) and incidence independent (for BCR) decompositions of the systems. Concentration robustness has also been studied in poly-PL and Hill-type systems \cite{FTJM2020, HEME2020}.

Our approach is based on a necessary and sufficient condition for ACR in PLP systems, and its analogue for BCR in CLP systems, in terms of their LP set's parameter subspaces' containment in species hyperplanes, i.e., subspaces of the form $\{ x \in \mathbb{R}^\mathscr{S} \mid x_X = 0 \}$ for a species $X$. A novel result is a necessary condition yielding an upper bound for the number of species in which the system admits concentration robustness. It also provides a simple procedure for determining species with concentration robustness. Since most of the systems mentioned above are LP systems, this approach provides a uniform view for many known results as well.

Our results on LP systems enable their novel use as ``control components'' for concentration robustness in decompositions (see Section \ref{sec:control}). LP systems of Shinar-Feinberg type, i.e., those with SF-pairs, led us to broaden the ``building blocks'' framework for constructing larger systems with concentration robustness studied in \cite{FOMF2021, FOME2021}. Besides the low deficiency building blocks already considered, higher deficiency weakly reversible LP systems with SF-pairs could be utilized in constructing appropriate decompositions. New computational approaches for determining independence presented by B. Hernandez and R. De la Cruz \cite{HEDC2021} and for detecting incidence independence by L. Fontanil and E. Mendoza \cite{FONME2022}, which will be useful for the broader framework, are also briefly reviewed.

In addition, in Section \ref{subsec:computational}, we state an ERRATUM for Theorem 2 of \cite{FOME2021}. It should be replaced by the said paper's more restrictive Theorem 6 which is independently proven in its Appendix.

The paper is organized as follows: Section \ref{sec:fundamentals} collects basic concepts and results on chemical reaction networks and kinetic systems needed in the later sections. In Section \ref{sec:criterion}, after the fundamentals of LP sets and LP systems are introduced, results on LP mass action systems and power law systems are reviewed before the necessary and sufficient condition (the species hyperplane criterion) for concentration robustness is derived. Section \ref{sec:control} discusses the role of LP subnetworks as ``control components'' in decompositions. In Section \ref{sec:block}, some benefits from adding LP systems of Shinar-Feinberg type to the ``building block'' framework for concentration robustness are presented. Section \ref{sec:concenRobust} provides some new computational approaches based on the results of the previous sections. Finally, summary and conclusion constitute Section \ref{sec:summary}.

Abbreviations used in this paper are listed in Table \ref{tab1} of the Appendix.


\section{Fundamentals of Chemical Reaction Networks and Kinetic Systems}
\label{sec:fundamentals}

In this section, we discuss fundamental concepts and results about chemical reaction networks and chemical kinetic systems. Moreover, we explore reaction networks as a digraph with vertex labeling.

\begin{definition}
	A \textbf{chemical reaction network} (CRN) is a digraph $(\mathscr{C}, \mathscr{R})$ where each vertex has positive degree and stoichiometry, i.e., there is a finite set $\mathscr{S}$ (whose elements are called \textbf{species}) such that $\mathscr{C}$ is a subset of $\mathbb{R}^\mathscr{S}_\geq$. Each vertex is called a \textbf{complex} and its coordinates in $\mathbb{R}^\mathscr{S}_\geq$ are called \textbf{stoichiometric coefficients}. The arcs are called \textbf{reactions}. We denote a CRN $\mathscr{N}$ as $\mathscr{N} = (\mathscr{S}, \mathscr{C}, \mathscr{R})$.
\end{definition}

In this paper, $\mathbb{R}^\mathscr{S}$, $\mathbb{R}^\mathscr{C}$, and $\mathbb{R}^\mathscr{R}$ denote the vector spaces of real-valued functions with domain $\mathscr{S}$, $\mathscr{C}$, and $\mathscr{R}$, respectively. A subscript of $\geq$ or $>$ denotes restriction to nonnegative or positive real numbers, respectively. Similarly, the sets of nonnegative and positive real numbers are denoted $\mathbb{R}_\geq$ and $\mathbb{R}_>$, respectively.

We denote the number of species with $m$, the number of complexes with $n$, and the number of reactions with $r$. We implicitly assume the elements of the sets are numbered and let $\mathscr{S} = \{X_1, \ldots, X_m\}$, $\mathscr{C} = \{C_1, \ldots, C_n\}$, and $\mathscr{R} = \{R_1, \ldots, R_r\}$, representing the set of species, complexes, and reactions, respectively. Furthermore, we denote the set of reactions as $\mathscr{R} \subset \mathscr{C} \times \mathscr{C}$. $(C_i, C_j) \in \mathscr{R}$ corresponds to the familiar notation $C_i \rightarrow C_j$.

Consider the reaction
\begin{equation*}
	\alpha X_1 + \beta X_2 \rightarrow \gamma X_3
\end{equation*}
where $X_1$, $X_2$, and $X_3$ are the species. The complexes are $\alpha X_1 + \beta X_2$ and $\gamma X_3$. In particular, $\alpha X_1 + \beta X_2$ is called the \textbf{reactant} (or \textbf{source}) \textbf{complex} and $\gamma X_3$ is the \textbf{product complex}. The number of reactant complexes is denoted by $n_r$. The stoichiometric coefficients are the nonnegative integer coefficients $\alpha$, $\beta$, and $\gamma$.

\begin{definition}
	The \textbf{reactant map} $\rho: \mathscr{R} \rightarrow \mathscr{C}$ maps a reaction $(C_i, C_j) \in \mathscr{R}$ to its reactant complex $C_i \in \mathscr{C}$.
\end{definition}

\begin{definition}
	Given a CRN $\mathscr{N} = (\mathscr{S}, \mathscr{C}, \mathscr{R})$, the \textbf{incidence map} $I_a: \mathbb{R}^\mathscr{R} \rightarrow \mathbb{R}^\mathscr{C}$ is a linear map such that for each reaction $R = (C_i, C_j) \in \mathscr{R}$, the basis vector $\omega_R$ is mapped to the vector $\omega_{C_j} - \omega_{C_i} \in \mathscr{C}$.
\end{definition}

\begin{definition}
	Given a CRN $(\mathscr{S}, \mathscr{C}, \mathscr{R})$, the \textbf{reactions map} $\rho': \mathbb{R}^\mathscr{C} \rightarrow \mathbb{R}^\mathscr{R}$ is given by $f: \mathscr{C} \rightarrow \mathbb{R}$ mapped to $f \circ \rho$ where $\rho$ is the reactant map.
\end{definition}

\begin{definition}
	The \textbf{stoichiometric subspace} of a reaction network $(\mathscr{S}, \mathscr{C}, \mathscr{R})$ is the linear subspace of $\mathbb{R}^\mathscr{S}$ given by $S = \text{span} \{C_j - C_i \in \mathbb{R}^\mathscr{S} \mid (C_i, C_j) \in \mathscr{R}\}$. The \textbf{rank} of the network is given by $s = \text{dim } S$. 
\end{definition}


\begin{definition}
	The \textbf{linkage classes} of a CRN are the subnetworks of its reaction graph where for any complexes $C_i$ and $C_j$ of the subnetwork, there is a path between them. The number of linkage classes is denoted by $\ell$.
\end{definition}

The linkage class is said to be a \textbf{strong linkage class} if there is a directed path from $C_i$ to $C_j$, and vice versa, for any complexes $C_i$ and $C_j$ of the subnetwork. The number of strong linkage classes is denoted by $s \ell$. Moreover, \textbf{terminal strong linkage classes}, the number of which is denoted as $t$, are the maximal strongly connected subnetworks where there are no edges (reactions) from a complex in the subgraph to a complex outside the subnetwork. Complexes belonging to terminal strong linkage classes are called \textbf{terminal}; otherwise, they are called \textbf{nonterminal}. 

\begin{example}
    \label{example1}
	Consider the following CRN:
	\begin{eqnarray'}
		&& R_1: 2 A_1 \rightarrow A_3	\\
		&& R_2: A_2 + A_3 \rightarrow A_3 \\
		&& R_3: A_3 \rightarrow A_2 + A_3 \\
		&& R_4: 3 A_4 \rightarrow A_2 + A_3 \\
		&& R_5: 2 A_1 \rightarrow 3 A_4.
	\end{eqnarray'}
		\indent We have
	\begin{eqnarray'}
		&& \mathscr{S} = \{A_1,A_2, A_3, A_4\} \\
		&& \mathscr{C} = \{C_1 = 2A_1, C_2 = A_2 + A_3, C_3 = A_3, C_4 = 3 A_4\}.
	\end{eqnarray'}
	Thus, there are $m = 4$ species, $n = 4$ complexes, $n_r = 4$ reactant complexes, and $r = 5$ reactions. The number of linkage classes is $\ell = 1$: $\{2 A_1, A_3, A_2 + A_3, 3 A_4\}$; the number of strong linkage classes is $s \ell = 3$: $\{A_3, A_2 + A_3\}, \{2 A_1\}, \{3 A_4\}$; and the number of terminal strong linkage classes is $t = 1$: $\{A_3, A_2 + A_3\}$. The rank of the CRN is $s = 3$.
\end{example}

\begin{definition}
	A CRN with $n$ complexes, $n_r$ reactant complexes, $\ell$ linkage classes, $s\ell$ strong linkage classes, and $t$ terminal strong linkage classes is called
	\begin{enumerate}[(i)]
		\item \textbf{weakly reversible} if $s\ell = \ell$;
		\item \textbf{$t$-minimal} if $t = \ell$;
		\item \textbf{point terminal} if $t = n - n_r$; and
		\item \textbf{cycle terminal} if $n - n_r = 0$.
	\end{enumerate}
\end{definition}

As observed in Example \ref{example1}, since $s\ell = 3 \neq 1 = \ell$, the network is not weakly reversible. $t = 1 = \ell$ implies that the network is $t$-minimal. Lastly, $t = 1 \neq 0 = 4 - 4 = n - n_r$ implies that the network is not point terminal but is cycle terminal.

\begin{definition}
	The \textbf{deficiency} of a CRN is the integer $\delta = n - \ell - s$ where $n$ is the number of complexes, $\ell$ is the number of linkage classes, and $s$ is the rank.
\end{definition}

In Example \ref{example1}, the deficiency of the network is $\delta = n - \ell - s = 4 - 1 - 3 = 0$.

\begin{definition}
	A \textbf{kinetics} of a CRN $\mathscr{N} = (\mathscr{S}, \mathscr{C}, \mathscr{R})$ is an assignment of a rate function $K_{C_i \rightarrow C_j}: \Omega_K \rightarrow \mathbb{R}_\geq$ to each reaction $(C_i, C_j) \in \mathscr{R}$ where $\Omega_K$ is a set such that $\mathbb{R}^\mathscr{S}_> \subseteq \Omega_K \subseteq \mathbb{R}^\mathscr{S}_\geq$ and
	\begin{equation*}
		K_{C_i \rightarrow C_j} (c) > 0 \text{ for all } c \in \Omega_K.
	\end{equation*}
	The kinetics of a network is denoted by $K = [K_1, \ldots, K_r]^T$. The pair $(\mathscr{N}, K)$ is called a \textbf{chemical kinetic system} (CKS).
\end{definition}

A kinetics gives rise to two closely related objects: the species formation rate function and the associated ordinary differential equation system.

\begin{definition}
	The \textbf{species formation rate function} (SFRF) of a CKS is defined as
	\begin{equation*}
		f(x) = N K(x) = \sum_{C_i \rightarrow C_j} K_{C_i \rightarrow C_j} (x) (C_j - C_i)
	\end{equation*}
	where $N$ is called the \textbf{stoichiometric matrix} and $K(x)$ is called the \textbf{kinetic vector} (or \textbf{kinetics}) of the CKS. The equation $\dot{x} = f(x)$ is the \textbf{ordinary differential equation} (ODE) system or \textbf{dynamical system} of the CKS.
\end{definition}

The dynamical system of the CRN in Example \ref{example1} can be written as
\begin{equation*}
	\dot{X} =
	\left[
		\begin{array}{c}
			\dot{A}_1 \\
			\dot{A}_2 \\
			\dot{A}_3 \\
			\dot{A}_4 \\
		\end{array}
	\right]
	=
	\left[
		\begin{array}{rrrrr}
		-2 & 0 & 0 & 0 & -2 \\
		0 & -1 & 1 & 1 & 0 \\
		1 & 0 & 0 & 1 & 0 \\
		0 & 0 & 0 & -3 & 3 \\
		\end{array}
	\right]
	\left[
		\begin{array}{l}
			k_1 A_1^{f_{11}} \\
			k_2 A_2^{f_{22}} A_3^{f_{23}} \\
			k_3 A_3^{f_{33}} \\
			k_4 A_4^{f_{44}} \\
			k_5 A_1^{f_{51}} \\
		\end{array}
	\right]
	= N K(x).
\end{equation*}
\begin{definition}
	The \textbf{set of positive equilibria} of a CKS $(\mathscr{N}, K)$ is given by
	\begin{equation*}
		E_+ (\mathscr{N}, K) = \{ x \in \mathbb{R}^\mathscr{S}_> \mid f(x) = 0 \}
	\end{equation*}
	where $f$ is the SFRF  of the CKS.
\end{definition}

Analogously, the \textbf{set of complex balanced equilibria} of a CKS $(\mathscr{N}, K)$ \cite{HOJA1972} is given by
\begin{equation*}
	Z_+ (\mathscr{N}, K) = \{ x \in \mathbb{R}^\mathscr{S}_> \mid I_a K(x) = 0 \} \subseteq E_+ (\mathscr{N}, K)
\end{equation*}
where $I_a$ is the incidence map.

A positive vector $c \in \mathbb{R}^\mathscr{S}$ is \textbf{complex balanced} if $K(c)$ is contained in the kernel $\text{Ker } I_a$, and a CKS is \textbf{complex balanced} if it has a complex balanced equilibrium.

An ODE system is under \textbf{power law kinetics} (PLK) if its kinetics has the form
\begin{equation*}
	K_i (x) = k_i \prod_{j=1}^m x_j^{f_{ij}} \text{ for } 1 \leq i \leq r
\end{equation*}
with $k_i \in \mathbb{R}_>$ (called the \textbf{rate constant} associated with reaction $R_i$) and $f_{ij} \in \mathbb{R}$ (called the \textbf{kinetic order} of species $x_j$). PLK is defined by an $r \times m$ matrix $F = [f_{ij}]$ called the \textbf{kinetic order matrix} and vector $k \in \mathbb{R}^r$ called the \textbf{rate vector}. We refer to a CRN with PLK as a \textbf{power law system}. A particular example of PLK is the well-known \textbf{mass action kinetics} (MAK) where the kinetic order matrix consists of stoichiometric coefficients of the reactants. We refer to a CRN with MAK as a \textbf{mass action system}.

In Example \ref{example1}, we assumed PLK so that the kinetic order matrix is
\begin{equation*}
	F =
	\left[
		\begin{array}{cccc}
			f_{11} & 0 & 0 & 0 \\
			0 & f_{22} & f_{23} & 0 \\
			0 & 0 & f_{33} & 0 \\
			0 & 0 & 0 & f_{44} \\
			f_{51} & 0 & 0 & 0 \\
		\end{array}
	\right].
\end{equation*}
	\indent We associate three linear maps to a positive element $k$ of $\mathbb{R}^\mathscr{R}$.

\begin{definition}
	For a reaction $R \in \mathscr{R}$, the \textbf{$k$-diagonal map} $\text{diag}(k)$ maps the basis vector $\omega_R$ to $k_R \omega_R$. The \textbf{$k$-incidence map} $I_k$ is defined as the composition $\text{diag}(k) \circ \rho'$ where $\rho'$ is the reactions map. The \textbf{$k$-Laplacian map} $A_k: \mathbb{R}^\mathscr{C} \rightarrow \mathbb{R}^\mathscr{C}$ is defined as the composition $A_k = I_a \circ I_k$ where $I_a$ is the incidence map.
\end{definition}

\begin{definition}
	A power law system has \textbf{reactant-determined kinetics} (of type \linebreak PL-RDK) if for any two reactions $R_i$ and $R_ j$ with identical reactant complexes, the corresponding rows of kinetic orders in $F$ are identical, i.e., $f_{ik} = f_{jk}$ for $k = 1, \ldots, m$. A power law system has \textbf{non-reactant-determined kinetics} (of type PL-NDK) if there exist two reactions with the same reactant complexes whose corresponding rows in $F$ are not identical.
\end{definition}


\section{A Species Hyperplane Criterion for Concentration Robustness in Systems of LP Type}
\label{sec:criterion}

In this section, we review the concept of a log-parametrized (LP) set in the species space $\mathbb{R}^\mathscr{S}$ of a network, and recall some results of Feinberg, Horn, and Jackson about such sets, which determine important properties of kinetic systems of LP type. We then provide an overview of the known examples of such systems. Finally, after reviewing the fundamentals of concentration robustness, we present simple criteria for the property in systems of LP type and illustrate it with some examples.

\subsection{Fundamentals of LP Sets and Systems of LP Type}
\label{subsec:fundamentals}

We begin with the concept of an LP set of a CKS.

\begin{definition}
	Given a set of species $\mathscr{S}$, a \textbf{log-parametrized} (LP) set is a non-empty subset of $\mathbb{R}^\mathscr{S}_>$ of the form $E(P, x^*) := \{x \in \mathbb{R}^\mathscr{S}_> \mid \log x - \log x^* \in P^\perp\}$ where $P$ (called the LP set's \textbf{flux subspace}) is a subspace of $\mathbb{R}^\mathscr{S}$, $x^*$ (called the LP set's \textbf{reference point}) is a given element of $\mathbb{R}^\mathscr{S}_>$, and $P^\perp$ (called the LP set's \textbf{parameter subspace}) is the orthogonal complement of $P$. The positive cosets of $P$ are called the LP set's \textbf{flux classes}.
\end{definition}

In \cite{FEIN1979}, M. Feinberg derived the following important property of an LP set based on the work by F. Horn and R. Jackson:

\begin{proposition}
	\label{prop:3.2}
	For an LP set $E = E(P, x^*)$ and any of its flux classes $Q$, $\vert E \cap Q \vert = 1$.
\end{proposition}

For a proof of Proposition \ref{prop:3.2}, see \cite{JMT2021}.

We now introduce concepts relating LP sets and equilibria sets of kinetic systems.

\begin{definition}
	A subset $E$ of the set of positive equilibria $E_+ (\mathscr{N}, K)$ of a CKS $(\mathscr{N}, K)$ is of \textbf{LP type} (or simply an \textbf{LP system}) if $E$ is an LP set, i.e., $E = E(P, x^*)$ for a subspace $P$ of $\mathbb{R}^\mathscr{S}$ and an element $x^* \in E$. A CKS is a \textbf{positive equilibria LP system} (of PLP type or simply a \textbf{PLP system}) if $E_+ (\mathscr{N}, K) \neq \varnothing$ and of LP type for a subspace $P_E$ of $\mathbb{R}^\mathscr{S}$. A CKS is a \textbf{complex balanced equilibria LP system} (of CLP type or simply a \textbf{CLP system}) if its set of complex balanced equilibria $Z_+ (\mathscr{N}, K) \neq \varnothing$ and of LP type for a subspace $P_Z$ of $\mathbb{R}^\mathscr{S}$. An LP system is a PLP or CLP system. It is a \textbf{bi-LP system} if it is both PLP and CLP with $P_E = P_Z$. $P_E$ and $P_Z$ are the LP system's flux subspaces, $P_E^\perp$ and $P_Z^\perp$ its parameter subspaces, and the positive cosets of $P_E$ and $P_Z$ are the LP system's flux classes.
\end{definition}

The following proposition justifies the term ``parameter subspace'' for $P_E^\perp$ and $P_Z^\perp$:

\begin{proposition}
	\label{prop:3.4}
	Let $(\mathscr{N}, K)$ be a CKS.
	\begin{enumerate}[(i)]
		\item If $(\mathscr{N}, K)$ is a PLP system with flux subspace $P_E$ and reference point $x^* \in E_+ (\mathscr{N}, K)$, then the map $L_{x^*}: E_+ (\mathscr{N}, K) \rightarrow P_E^\perp$ given by $L_{x^*} (x) = \log x - \log x^*$ is a bijection.
		\item If $(\mathscr{N}, K)$ is a CLP system with flux subspace $P_Z$ and reference point $x^* \in Z_+ (\mathscr{N}, K)$, then the restriction to $Z_+ (\mathscr{N}, K)$ of $L_{x^*}: Z_+ (\mathscr{N}, K) \rightarrow P_Z^\perp$ is a bijection.
	\end{enumerate}
\end{proposition}

For a proof of Proposition \ref{prop:3.4}, see \cite{JMT2021}.

The following theorem collects the important properties of LP systems derived from the results of Feinberg, Horn, and Jackson reviewed above:

\begin{theorem}
	\label{thm:3.5}
	Let $(\mathscr{N}, K)$ be a CKS.
	\begin{enumerate}[(i)]
		\item If $(\mathscr{N}, K)$ is a PLP system, then $\vert E_+ (\mathscr{N}, K) \cap Q \vert = 1$ for any of its flux classes $Q$.
		\item If $(\mathscr{N}, K)$ is a CLP system, then $\vert Z_+ (\mathscr{N}, K) \cap Q \vert = 1$ for any of its flux classes $Q$.
		\item If $(\mathscr{N}, K)$ is a bi-LP system, then it is \textbf{absolutely complex balanced}, i.e., $E_+ (\mathscr{N}, K) = Z_+ (\mathscr{N}, K)$.
	\end{enumerate}
\end{theorem}

For a proof of Theorem \ref{thm:3.5}, see \cite{JMT2021}.

The first (and best-known) example of an LP system is the set of complex balanced mass action systems studied by F. Horn and R. Jackson in 1972 \cite{HOJA1972}. In this case, as shown in the succeeding formulation of their result in M. Feinberg's 1979 lectures \cite{FEIN1979}, the LP system property is equivalent to other significant characteristics of the system.

\begin{theorem}
	Let $\mathscr{N} = (\mathscr{S}, \mathscr{C}, \mathscr{R})$ be a CRN. Suppose $(\mathscr{N}, K)$ is a mass action system with stoichiometric subspace $S$ and SFRF $f$. If there exists $c^* \in \mathbb{R}^\mathscr{S}_>$ such that
	\begin{equation*}
		A_k \psi_K (c^*) = 0,
	\end{equation*}
	where $A_k$ is the $k$-Laplacian map and $\psi_K$ is a factor map of $K$, then
	\begin{equation*}
		f(c) (\log c - \log c^*) \leq 0 \text{ for all } c \in \mathbb{R}^\mathscr{S}_>.
	\end{equation*}
	Moreover, for $c \in \mathbb{R}^\mathscr{S}_>$, the following are equivalent:
	\begin{enumerate}[(i)]
		\item $f(c) (\log c - \log c^*) = 0$;
		\item $\log c - \log c^* \in S^\perp$;
		\item $A_k \psi_K (c) = 0$; and
		\item $f(c) = 0$.
	\end{enumerate}
\end{theorem}

Statements $(ii)$ and $(iii)$ establish that the set of complex balanced mass action systems are CLP systems with $P_Z = S$. Furthermore, the inference $(iv) \Rightarrow (iii)$ shows that the systems are absolutely complex balanced, i.e., every positive equilibrium is complex balanced. Hence, they are also bi-LP since $P_E = P_Z$.

The other class of mass action systems with strong LP characteristics are those with independent linkage classes. This is equivalent to $\delta = \delta_1 + \ldots + \delta_\ell$ where $\delta$ is the deficiency of the network and $\delta_i$ is the deficiency of linkage class $i$ \cite{FML2021}. B. Boros \cite{BOROS2013} lists the following known relationship between $E_+ (\mathscr{N}, K)$ and $E(S, x^*)$ in general:

\begin{proposition}
	\label{prop:3.7}
	Let $(\mathscr{N}, K)$ be a mass action system with stoichiometric subspace $S$, deficiency $\delta = \delta_1 + \ldots + \delta_\ell$ where $\delta_i$ is the deficiency of linkage class $i$, and set of positive equilibria $E_+ (\mathscr{N}, K) \neq \varnothing$. Then the LP set $E(S, x^*) \subseteq E_+ (\mathscr{N}, K)$ for any $x^* \in E_+ (\mathscr{N}, K)$.
\end{proposition}

According to the Deficiency One Theorem which M. Feinberg proved in 1995 \cite{FEIN1995}, if for a weakly reversible mass action system with deficiency $\delta = \delta_1 + \ldots + \delta_\ell$ (where $\delta_i$ is the deficiency of linkage class $i$) we have $\delta_i \leq 1$ for all $i$, then equality holds in Proposition \ref{prop:3.7}, i.e., $E(S, x^*) = E_+ (\mathscr{N}, K)$ so that $(\mathscr{N}, K)$ is a PLP system. B. Boros extended this result to a class of $t$-minimal networks by providing a necessary and sufficient condition for such systems to have a positive equilibrium \cite{BOROS2012}.

\subsection{A Review of Power Law Systems of LP Type}
\label{subsec:reviewPowerLaw}

In this section, we collect the known results on LP systems among PL-RDK and PL-NDK systems. In 2014, S. M\"{u}ller and G. Regensburger introduced generalized mass action systems in \cite{MURE2012} as follows:

\subsubsection{PL-RDK Systems of LP Type}
\label{subsubsec:plrdk}

\begin{definition}
	A \textbf{generalized mass action system} (GMAS) is a triple $(G, \varphi, \tilde{\varphi})$ where $G = (V, E)$ is a digraph ($V$ the set of vertices and $E$ the set of edges), $\varphi: V \rightarrow \mathbb{R}^m$, and $\tilde{\varphi}: V_S \rightarrow \mathbb{R}^m$ ($V_S$ is the subset of source vertices and $m$ is the number of species in the system).
\end{definition}

The set of PL-RDK systems can be mapped bijectively to the subset of GMAS where $\varphi$ is injective and $\tilde{\varphi}$ maps the zero complex (if it is a source) to the zero vector in $\mathbb{R}^m$. In this bijection, the subset of factor span surjective systems is mapped to the GMAS where, additionally, $\tilde{\varphi}$ is injective.

One of the key concepts of the GMAS theory is that of the \textbf{kinetic order subspace} $\tilde{S} := \text{span} \{ \tilde{\varphi}(v') - \tilde{\varphi}(v) \}$ for any edge $v \rightarrow v'$ of a cycle terminal digraph, i.e., $V_S = V$. S. M\"{u}ller and G.  Regensburger showed that, analogous to mass action systems, any complex balanced PL-RDK system is a CLP system with $P_Z = \tilde{S}$. However, in contrast to mass action systems, only special subsets of PL-RDK systems with positive deficiency are absolutely complex balanced. They also showed that, for a weakly reversible \linebreak PL-RDK system, there are rate constants for which the system is complex balanced. In particular, it is complex balanced for all rate constants if and only if its \textbf{kinetic deficiency} $\tilde{\delta} := n - \ell - \text{dim } \tilde{S}$ is 0 ($n$ is the number of complexes and $\ell$ is the number of linkage classes in the system).

For PL-RDK systems with independent linkage classes (or, equivalently, with deficiency $\delta = \delta_1 + \ldots + \delta_\ell$ where $\delta_i$ is the deficiency of linkage class $i$), Talabis et al \cite{TAM2018} identified a subset of PL-RDK systems called \textbf{$\hat{T}$-rank maximal kinetic} (PL-TIK) \textbf{systems}, for which the full analogue of the Deficiency One Theorem for mass action systems, including Boros's criterion, is valid. This class of power law systems is defined with the help of the $T$-matrix, a $\vert V_S \vert \times m$ matrix whose columns are the images $\tilde{\varphi}(v)$. A kinetics is PL-TIK if and only if the augmented $T$-matrix $\hat{T}$, constructed by adding rows of the characteristic functions of the linkage classes, has maximal column rank (see \cite{TMJ2019} for details). In particular, any PL-TIK on a $t$-minimal network with independent linkage classes of low deficiency ($\delta = 0$ or $\delta = 1$) satisfying Boros's criterion is PLP with $P_E = \tilde{S}_R$ where $\tilde{S}_R$ is called the \textbf{kinetic reactant flux subspace}. The latter is defined as the image of $T I_{a,R}$ where $T$ is the $T$-matrix and $I_{a,R}$ is the restriction of the incidence matrix to the maximal cycle terminal subnetwork. If the network is cycle terminal, then $\tilde{S}_R = \tilde{S}$ follows from this and from \cite{TMJ2019} that any weakly reversible PL-TIK system satisfying the Deficiency One Theorem for PL-TIK is bi-PL with LP flux subspace $\tilde{S}$.

\subsubsection{PL-NDK Systems of LP Type}
\label{subsubsec:plndk}

Recall that for a CRN $\mathscr{N} = (\mathscr{S}, \mathscr{C}, \mathscr{R})$ a \textbf{covering} of $\mathscr{N}$ is a collection of subsets $\{\mathscr{R}_1, \ldots, \mathscr{R}_k\}$ whose union is $\mathscr{R}$. A covering is called a \textbf{decomposition} of $\mathscr{N}$ if the sets $\mathscr{R}_i$ form a partition of $\mathscr{R}$. Alternatively, $\mathscr{N}$ has a decomposition into subnetworks $\mathscr{N}_i = (\mathscr{S}_i, \mathscr{C}_i, \mathscr{R}_i)$ for $i = 1, \ldots, k$ if $\displaystyle \mathscr{S} = \bigcup_{i=1}^k \mathscr{S}_i$, $\displaystyle \mathscr{C} = \bigcup_{i=1}^k \mathscr{C}_i$, $\displaystyle \mathscr{R} = \bigcup_{i=1}^k \mathscr{R}_i$, and $\mathscr{R}_i \cap \mathscr{R}_j = \varnothing$ for $i \neq j$. We denote a decomposition of $\mathscr{N}$ into $k$ subnetworks as $\mathscr{N} = \mathscr{N}_1 \cup \ldots \cup \mathscr{N}_k$.

For any decomposition, the stoichiometric subspace of $\mathscr{N}$ is $\displaystyle S = \sum_{i=1}^k S_i$ where $S_i$ is the stoichiometric subspace of subnetwork $\mathscr{N}_i$, the image of the incidence map is \linebreak $\displaystyle \text{Im } I_a = \sum_{i=1}^k \text{Im } I_{a,i}$ where $I_{a,i}$ is the incidence map of subnetwork $\mathscr{N}_i$, and the decomposition is independent (incidence independent) if the first (second) sum is direct. A decomposition has a network property (e.g., weakly reversible) or kinetic property (e.g., PL-RDK) if each subnetwork has the said property.

In \cite{HEME2021}, PL-NDK systems with weakly reversible PL-RDK decompositions were investigated and the following application was derived:

\begin{theorem}
	Let $(\mathscr{N}, K)$ be a power law system with a weakly reversible PL-RDK decomposition $\mathscr{D}: \mathscr{N} = \mathscr{N}_1 \cup \ldots \cup \mathscr{N}_k$. If $\mathscr{D}$ is bi-level independent and of PLP type with $P_{E,i} = \tilde{S}_i$ where $P_{E,i}$ is the flux subspace of subnetwork $\mathscr{N}_i$ and $\tilde{S}_i$ is the kinetic order subspace of subnetwork $\mathscr{N}_i$, then $(\mathscr{N}, K)$ is a weakly reversible PLP system with flux subspace $\displaystyle P_E = \sum_{i=1}^k \tilde{S}_i$.
\end{theorem}

\textbf{Bi-level independent} means that both the decomposition as well as the induced decomposition of kinetic complexes (the kinetic complexes of a PL-RDK system are the images of the reactant complexes under the kinetic map $\tilde{\varphi}$ of S. M\"{u}ller and G. Regensburger) are independent.

This result is a broad generalization of the Deficiency Zero Theorem of Fortun et al \cite{FLRM2019}. The analogous result for complex balanced PL-NDK systems is contained in the following result:

\begin{theorem}
	Let $(\mathscr{N}, K)$ be a weakly reversible power law system with a complex balanced PL-RDK decomposition $\mathscr{D}: \mathscr{N} = \mathscr{N}_1 \cup \ldots \cup \mathscr{N}_k$ with $P_{Z,i} = \tilde{S}_i$ where where $P_{Z,i}$ is the flux subspace of subnetwork $\mathscr{N}_i$ and $\tilde{S}_i$ is the kinetic order subspace of subnetwork $\mathscr{N}_i$. If $\mathscr{D}$ is incidence independent and the induced covering is independent, then $(\mathscr{N}, K)$ is a weakly reversible CLP system with flux subspace $\displaystyle P_Z = \sum_{i=1}^k \tilde{S}_i$.
\end{theorem}

\subsection{The Species Hyperplane Containment Criterion for Concentration Robustness in LP Systems}
\label{subsec:hyperplane}

In this section, after a brief review of the fundamentals of concentration robustness, we derive a simple criterion, i.e., necessary and sufficient condition, for absolute concentration robustness (ACR) in PLP systems. The result is that the dimension of the LP set's flux subspace is an upper bound for the number of species with concentration robustness. The criterion also leads to simple computational algorithms for checking concentration robustness in LP systems. These are discussed in detail in Section \ref{sec:block}.

The concept of ACR was introduced by G. Shinar and M. Feinberg in \textit{Science} in 2010 \cite{SHFE2010} as follows: A CKS $(\mathscr{N}, K)$ has ACR in a species $X$ if the value of $X$ is identical for all positive equilibria of the system. Abstracting from earlier observations by biologists in systems in \textit{Escherichia coli}, they derived a sufficient condition for ACR in a species for deficiency one mass action systems (extensions of this result will be discussed in detail in Section \ref{sec:control}). The paper aroused interest in robustness research which included the work by Karp et al \cite{KMDDG2012} proposing further concentration robustness concepts such as bounded concentration robustness.

A different type of extension was introduced by N. Fortun and E. Mendoza in \cite{FOME2021}: a complex balanced CKS has balanced concentration robustness (BCR) in a species $X$ if the species has the same value at all complex balanced equilibria. In a general kinetic system, in contrast to a mass action system, there may be positive equilibria which are not complex balanced, so that BCR is a weaker property than ACR. We introduce a further generalization.

\begin{definition}
	A CKS is \textbf{concentration robust} for a subset $E$ of the set of positive equilibria in a species $X$ if $X$ has the same value for all equilibria in $X$. A general notation for this property is ``$(E, X)$ concentration robustness''.
\end{definition}

A CKS has ACR in a species $X$ if and only if all of its positive equilibria lie in the hyperplane $x_X = c$ (the so-called  ``ACR hyperplane''), where $c$ is a positive constant. This is equivalent to the differences of their logarithms lying in the species hyperplane $x_X = 0$,  due to the bijectivity of the logarithm on the positive real axis. Hence, a general ACR hyperplane criterion is that the $\text{span} \{ \log x - \log x^*  \mid x, x^* \in E_+ (\mathscr{N}, K) \}$, where $E_+ (\mathscr{N}, K)$ is the set of positive equilibria, lies in the species hyperplane $x_X = 0$. For LP sets, this span is equal to the parameter subspace and is easily computed in terms of network structures as shown below.

We now derive our \textbf{species hyperplane containment criterion} (or simply \textbf{species hyperplane criterion}) for concentration robustness.

\begin{theorem}
	\label{thm:hyperplanecriterion}
	Let $(\mathscr{N}, K)$ be an LP system.
	\begin{enumerate}[(i)]
		\item If $(\mathscr{N}, K)$ is a PLP system, then it has ACR in a species $X$ if and only if its parameter subspace $P_E^\perp$ is a subspace of the species hyperplane $\{ x \in \mathbb{R}^\mathscr{S} \mid x_X = 0 \}$. Similarly, if $(\mathscr{N}, K)$ is a weakly reversible CLP system, then it has BCR in a species $X$ if and only if its parameter subspace $P_Z^\perp$ is a subspace of the species hyperplane $\{ x \in \mathbb{R}^\mathscr{S} \mid x_X = 0 \}$.
		\item If $m_\text{ACR}$ is the number of species with ACR, then $m_\text{ACR} \leq \text{dim } P_E$ where $P_E$ is the flux subspace. If $m_\text{BCR}$ is the number of species with BCR, then $m_\text{BCR} \leq \text{dim } P_Z$ where $P_Z$ is the flux subspace.
	\end{enumerate}
\end{theorem}

\textit{Proof}.

$(i)$ (``$\Leftarrow$'') Since for any $p \in P_E^\perp$ we have $p_X = 0$, then for any $x \in E_+ (\mathscr{N}, K)$, $\log x_X = \log x_X^*$, and by the bijectivity of the logarithm $\log: \mathbb{R}^\mathscr{S}_> \rightarrow \mathbb{R}^\mathscr{S}$, we have $x_X = x_X^*$. Conversely (``$\Rightarrow$''), suppose the system has ACR in $X$. Then for any $x \in E_+ (\mathscr{N}, K)$, $x_X = x_X^*$, which implies that $\log x_X = \log x_X^*$. In view of the parametrization, this implies that $p_X = 0$ for any $p \in P_E^\perp$. An analogous argument clearly holds for CLP systems, complex balanced equilibria, and $P_Z^\perp$ elements.

$(ii)$ It follows from $(i)$ that if the system has ACR in the species $X_1, \ldots, X_{m_\text{ACR}}$, then $P_E^\perp$ is contained in $\displaystyle \bigcap_{i=1}^{m_{\text{ACR}}} \{ x \in \mathbb{R}^\mathscr{S} \mid x_{X_i} = 0 \}$. Hence, $\text{dim } P_E^\perp \leq m - m_\text{ACR}$ or $m_\text{ACR} \leq m - \text{dim } P_E^\perp = m - (m - \text{dim } P_E) = \text{dim } P_E$. An analogous inequality clearly holds for $m_\text{BCR}$. \hfill $\Box$

\begin{definition}
	A CRN $\mathscr{N} = (\mathscr{S}, \mathscr{C}, \mathscr{R})$ with stoichiometric subspace $S$ is said to be \textbf{conservative} if there exists an $m \in S^\perp$ such that $m \in \mathbb{R}^\mathscr{S}_>$ or, equivalently, $S^\perp \cap \mathbb{R}^\mathscr{S}_> \neq \varnothing$. Otherwise, it is said to be \textbf{nonconservative}.
\end{definition}

Our species hyperplane criterion has the following interesting consequences for mass action systems:

\begin{corollary}
	\textcolor{white}{}
	\begin{enumerate}[(i)]
		\item An LP mass action system with ACR or BCR in at least one species is nonconservative.
		\item A conservative mass action system with ACR or BCR in at least one species is not an LP system.
		\item A complex balanced conservative mass action system does not have ACR or BCR in any species.
	\end{enumerate}
\end{corollary}

\textit{Proof}.

$(i)$ follows from  the fact that the flux subspace is $S$, and hence every vector of the parameter subspace $S^\perp$ has at least one zero coordinate.

$(ii)$ is just a reformulation of $(i)$.

$(iii)$ follows from $(i)$ and the result of F. Horn and R. Jackson that any complex balanced mass action system is a CLP system. \hfill $\Box$


\section{LP Subsystems as Control Components for Concentration Robustness in Decompositions}
\label{sec:control}

In this section, we first use the species hyperplane criterion to construct a simple computational procedure for determining concentration robustness in a species of an LP system. We then discuss how the procedure enables the role of an LP subsystem as a ``control component'' for concentration robustness in a decomposition.

\subsection{A Simple Computational Approach to Concentration Robustness in LP Systems}
\label{subsec:computational}

The following proposition provides a simple computational procedure to determine for which species concentration robustness in an LP system holds:

\begin{proposition}
	\label{prop:4.1}
	Let $\{ v_1, \ldots, v_E \}$ be a basis of the parameter subspace $P_E^\perp$ of a PLP system $(\mathscr{N}, K)$. The system has ACR in species $X$ if and only if the coordinate corresponding to $X$ in each basis vector $v_{i,X} = 0$ for each $i = 1, \ldots, E$. Similarly, if $\{ w_1, \ldots, w_Z \}$ is a basis of the parameter subspace $P_Z^\perp$ of a CLP system, then the system has BCR in species $X$ if and only if the coordinate corresponding to $X$ in each basis vector $w_{i,X} = 0$ for each $i = 1, \ldots, Z$.
\end{proposition}

The proof is straightforward. Hence, one only needs to construct a basis for the parameter subspace $P_E^\perp$ or $P_Z^\perp$ to determine which species in the system has ACR or BCR.

\begin{example}
	For complex balanced mass action systems, the dimension of the flux subspace of its LP set $E$ is $\text{dim } P_E = s$ where $s$ is the rank of the system, so that the number of species with ACR is $m_\text{ACR} \leq s$.
\end{example}

This is an interesting relationship which, to our knowledge, has not been published to date. If a complex balanced mass action system has a unique positive equilibrium, then it has ACR in all species, so that $m_\text{ACR} = m$ (the number of species in the system). By Proposition \ref{prop:4.1}, it follows that $s = m$. Hence, the system is open, i.e., the linear space generated by the reaction vectors is the whole species space $\mathbb{R}^\mathscr{S}$. A conservative complex balanced mass action system, therefore, has at most $m - 1$ species with ACR.

\begin{example}
	\label{example2}
	Fortun et al \cite{FLRM2019} showed that the following subnetwork of Schmitz's global pre-industrial carbon cycle model is a PLP system with flux subspace $P_E = \tilde{S} = \tilde{S}_1 \oplus \tilde{S}_2$ (where $\tilde{S}$ and $\tilde{S}_i$ are the kinetic order subspaces of $\mathscr{N}$ and $\mathscr{N}_i$, respectively) defined by the following decomposition $\mathscr{D}$:
	\begin{equation*}
		\mathscr{N} = \mathscr{N}_1 \cup \mathscr{N}_2
	\end{equation*}
	where
	\begin{eqnarray'}
		&& \mathscr{N}_1 = \{ r_1: M_5 \rightarrow M_1, r_2: M_1 \rightarrow M_5, r_3: M_5 \rightarrow M_6, r_4: M_6 \rightarrow M_1 \} \\
		&& \mathscr{N}_2 = \{ r_5: M_2 \rightarrow M_1, r_6: M_4 \rightarrow M_2, r_7: M_1 \rightarrow M_3, r_8: M_3 \rightarrow M_4 \}.
	\end{eqnarray'}
	 The rank of $\mathscr{N}$, $\mathscr{N}_1$, and $\mathscr{N}_2$ are $s = 5$, $s_1 = 2$, and $s_2 = 3$, respectively. This implies that the decomposition is independent.
	 
	 The kinetic order subspaces of the subnetworks induced by $\mathscr{D}$ are
	\begin{eqnarray'}
		&& \tilde{\mathscr{N}}_{\mathscr{D},1} = \{ \tilde{r}_1: M_5 \rightarrow 0.36 M_1, \tilde{r}_2: 0.36 M_1 \rightarrow M_5, r_3: M_5 \rightarrow M_6, \tilde{r}_4: M_6 \rightarrow 0.36 M_1 \} \\
		&& \tilde{\mathscr{N}}_{\mathscr{D},2} = \{ \tilde{r}_5: 9.4 M_2 \rightarrow M_1, \tilde{r}_6: M_4 \rightarrow 9.4 M_2, r_7: M_1 \rightarrow M_3, r_8: M_3 \rightarrow M_4 \}.
	\end{eqnarray'}
	Set $\tilde{\mathscr{N}}_\mathscr{D} = \tilde{\mathscr{N}}_{\mathscr{D},1} \cup \tilde{\mathscr{N}}_{\mathscr{D},2}$. The number of complexes is $\tilde{n} = 7$ and the number of linkage classes is $\tilde{\ell} = 2$. So, $\tilde{n} - \tilde{\ell} = 5$. The rank of $\tilde{\mathscr{N}}_\mathscr{D}$, $\tilde{\mathscr{N}}_{\mathscr{D},1}$, and $\tilde{\mathscr{N}}_{\mathscr{D},2}$ are $\tilde{s} = 5$, $\tilde{s}_1 = 2$, and $\tilde{s}_2 = 3$, respectively. This implies that $\tilde{\mathscr{N}}_\mathscr{D} = \tilde{\mathscr{N}}_{\mathscr{D},1} \cup \tilde{\mathscr{N}}_{\mathscr{D},2}$ is an independent decomposition.
	
	We now apply Proposition \ref{prop:4.1}. The flux subspace induced by $\mathscr{D}$ is $\tilde{S}$ and we wish to determine a basis for $\tilde{S}^\perp$ to identify which species has ACR.
	
	Note, first of all, that the set of species of $\tilde{\mathscr{N}}_\mathscr{D}$ is $\mathscr{S} = \{ M_1, M_2, M_3, M_4, M_5, M_6 \}$.
	
	$\tilde{S}$ is represented by the following stoichiometric matrix:
	\begin{equation*}
		\tilde{N}_\mathscr{D} =
		\left[
		\begin{array}{rrrrrrrr}
			0.36 & -0.36 & 0 & 0.36 & 1 & 0 & -1 & 0 \\
			0 & 0 & 0 & 0 & -9.4 & 9.4 & 0 & 0 \\
			0 & 0 & 0 & 0 & 0 & 0 & 1 & -1 \\
			0 & 0 & 0 & 0 & 0 & -1 & 0 & 1 \\
			-1 & 1 & -1 & 0 & 0 & 0 & 0 & 0 \\
			0 & 0 & 1 & -1 & 0 & 0 & 0 & 0 \\
		\end{array}
		\right].
	\end{equation*}
	The columns of $\tilde{N}_\mathscr{D}$ are vectors in $\mathbb{R}^\mathscr{S}$ (or $\mathbb{R}^6$, loosely speaking, keeping in mind the order in which we listed the species above). This makes $\tilde{S}$ a linear subspace of $\mathbb{R}^6$.
	
	We take the transpose of $\tilde{N}_\mathscr{D}$ to get the ``reaction matrix'' $\tilde{R}_\mathscr{D}$ where each row represents a reaction:
	\begin{equation*}
		\tilde{R}_\mathscr{D} = \tilde{N}_\mathscr{D}^T =
		\left[
		\begin{array}{rrrrrr}
			0.36 & 0 & 0 & 0 & -1 & 0 \\
			-0.36 & 0 & 0 & 0 & 1 & 0 \\
			0 & 0 & 0 & 0 & -1 & 1 \\
			0.36 & 0 & 0 & 0 & 0 & -1 \\
			1 & -9.4 & 0 & 0 & 0 & 0 \\
			0 & 9.4 & 0 & -1 & 0 & 0 \\
			-1 & 0 & 1 & 0 & 0 & 0 \\
			0 & 0 & -1 & 1 & 0 & 0 \\
		\end{array}
		\right].
	\end{equation*}
		\indent A basis $B$ for $\tilde{R}_\mathscr{D}$ consists of reactions $\tilde{r}_1, r_3, \tilde{r}_5, \tilde{r}_6, r_7$:
	\begin{equation*}
		B =
		\left[
		\begin{array}{rrrrrr}
		0.36 & 0 & 0 & 0 & -1 & 0 \\
		0 & 0 & 0 & 0 & -1 & 1 \\
		1 & -9.4 & 0 & 0 & 0 & 0 \\
		0 & 9.4 & 0 & -1 & 0 & 0 \\
		-1 & 0 & 1 & 0 & 0 & 0 \\
		\end{array}
		\right].
	\end{equation*}
		\indent We now want to determine the orthogonal complement $\tilde{S}^\perp = \{x \in \mathbb{R}^6 \mid Bx = 0\}$. To do this, we write the augmented matrix $[B \mid 0]$ in its reduced row echelon form (rref). In the following, we do not include the 0 column from the augmented matrix. The species corresponding to each column are also labeled:
	\begin{equation*}
		\text{rref } B =
		\begin{blockarray}{cccccc}
			\textcolor{white}{0}M_1 & M_2 & M_3 & M_4 & M_5 & M_6 \\
			\begin{block}{[cccccr]}
				\textcolor{white}{0}1 & 0 & 0 & 0 & 0 & -2.7778\textcolor{white}{0} \bigstrut[t] \\
				\textcolor{white}{0}0 & 1 & 0 & 0 & 0 & -0.2955\textcolor{white}{0} \\
				\textcolor{white}{0}0 & 0 & 1 & 0 & 0 & -2.7778\textcolor{white}{0} \\
				\textcolor{white}{0}0 & 0 & 0 & 1 & 0 & -2.7778\textcolor{white}{0} \\
				\textcolor{white}{0}0 & 0 & 0 & 0 & 1 & -1\textcolor{white}{0} \bigstrut[b] \\
			\end{block}
		\end{blockarray}.
	\end{equation*}
		\indent Using $B$ in rref, we write the matrix in its system of equations form. We also express the pivots in terms of the nonpivots (the column for $M_6$ is a nonpivot column):
	\begin{eqnarray*}
		M_1	&=& 2.7778 M_6 \\
		M_2	&=& 0.2955 M_6 \\
		M_3	&=& 2.7778 M_6 \\
		M_4	&=& 2.7778 M_6 \\
		M_5	&=& M_6.
	\end{eqnarray*}
	Therefore, the desired orthogonal complement is the set
	\begin{equation*}
		\tilde{S}^\perp =
		\left\{
			\left[
				\begin{array}{l}
					M_1 \\
					M_2 \\
					M_3 \\
					M_4 \\
					M_5 \\
					M_6 \\
				\end{array}
			\right]
			=
			\left[
				\begin{array}{r}
					2.7778 \\
					0.2955 \\
					2.7778 \\
					2.7778 \\
					1 \\
					1 \\
				\end{array}
			\right] M_6
		\right\}.
	\end{equation*}
	The vector beside the nonpivot is the basis vector for $\tilde{S}^\perp$. Since there is no zero coordinate for the vector, the species $M_1, \ldots, M_6$ have no ACR in the system.
	
\end{example}

\noindent
\textbf{ERRATUM.} The result in Example \ref{example2} contradicts the conclusion in Example 3 of \cite{FOME2021}, indicating an error in the scope of Theorem 2 of the paper. The theorem should be restricted to PL-RDK systems where the reactions of the Shinar-Feinberg pair lie in the same linkage class (as formulated in Theorem 6 which is independently proven in the Appendix of the paper). Consequently, the statements for deficiency zero building blocks in Propositions 7 and 8 in the paper need to be adjusted. Furthermore, the algorithms based on these propositions in \cite{FOMF2021} need to be refined to meet the restriction.

\subsection{The ``Control Component'' Role of LP Subsystems in a Decomposition for Concentration Robustness}
\label{subsec:control}

In the next proposition, we collect some important properties of concentration robustness for a species in a decomposition subnetwork, particularly in relation to its concentration robustness in the whole network.

\begin{proposition}
	\label{prop:4.4}
	Let $\mathscr{N} = \mathscr{N}_1 \cup \ldots \cup \mathscr{N}_k$ be a decomposition of the CRN $\mathscr{N} = (\mathscr{S}, \mathscr{C}, \mathscr{R})$ with $\mathscr{N}_i = (\mathscr{S}_i, \mathscr{C}_i, \mathscr{R}_i)$ for $i = 1, \ldots, k$. Denote by $\mathscr{S}_\text{ACR}$ ($\mathscr{S}_\text{BCR}$) and $\mathscr{S}_{\text{ACR},i}$ ($\mathscr{S}_{\text{BCR},i}$) the sets of species with ACR (BCR) of the network $\mathscr{N}$ and the subnetwork $\mathscr{N}_i$, respectively.
	\begin{enumerate}[(i)]
		\item If species $X$ has concentration robustness in $\mathscr{N}_i$, then $X \in \mathscr{S}_i$, i.e., $\vert \mathscr{S}_{\text{ACR},i} \vert \leq \vert \mathscr{S}_i \vert$ and $\vert \mathscr{S}_{\text{BCR},i} \vert \leq \vert \mathscr{S}_i \vert$.
		\item If species $X$ has ACR in $\mathscr{N}_i$ and the decomposition is independent, then $X$ has ACR in $\mathscr{N}$, i.e., $\displaystyle \vert \mathscr{S}_{\text{ACR},i} \vert \leq \Biggl| \bigcup_{i=1}^k \mathscr{S}_{\text{ACR},i} \Biggr| \leq \vert \mathscr{S}_\text{ACR} \vert$.
		\item If species $X$ has BCR in $\mathscr{N}_i$ and the decomposition is incidence independent, then $X$ has BCR in $\mathscr{N}$, i.e., $\displaystyle \vert \mathscr{S}_{\text{BCR},i} \vert \leq \Biggl| \bigcup_{i=1}^k \mathscr{S}_{\text{BCR},i} \Biggr| \leq \vert \mathscr{S}_\text{BCR} \vert$.
	\end{enumerate}
\end{proposition}

\textit{Proof}.

$(i)$ If a species does not occur in a subnetwork, since the set of equilibria for the subsystem includes elements with all positive values in that species, then concentration robustness in that species is impossible.

$(ii)$ If the decomposition is independent, then the theorem of M. Feinberg says that the set of equilibria of the whole system is the intersection of the equilibria sets of the subnetworks. Hence, the set of equilibria of the whole network is a subset of the set of equilibria for each subnetwork, implying the claim.

$(iii)$ The argument is analogous in view of the theorem of Farinas et al \cite{FML2020, FML2021} for complex balanced equilibria under incidence independent decompositions. \hfill $\Box$

If a subnetwork is an LP system, the simple procedure for determining which of the species occurring in that subsystem have concentration robustness allows the easy computation of the sets and bounds given in the proposition. This provides us a way of checking the ``amount'' of concentration robustness occurring in the system. The higher the presence of LP subnetworks in a decomposition, the more accurate the ``control'' of concentration robustness can be effected.

\begin{remark}
	In Example \ref{example2} and in the generalization of the Deficiency Zero Theorem of Fortun et al derived in \cite{FLRM2019}, all subnetworks are PLP systems, so that the entire network is also a PLP system. However, in cases where not all subnetworks are LP systems, the set of species in those which are LP systems still provides a useful lower bound for the set of species with ACR for the whole network.
\end{remark}


\section{Concentration Robustness in LP Power Law Systems of SF-type}
\label{sec:block}

In this section, we present a further application of the species hyperplane criterion for concentration robustness in LP systems by broadening the ``building blocks'' framework for power law systems introduced by N. Fortun and E. Mendoza in \cite{FOME2021}. Specifically, we add LP systems with Shinar-Feinberg reaction pairs (SF-pairs) of arbitrary deficiency to the original low deficiency ($\delta = 0$ or $\delta = 1$) building blocks in the framework. This addition will allow easier identification of concentration robustness in larger networks with arbitrary deficiency.

\subsection{A Review of Concentration Robustness in Power Law Systems of SF-Type}
\label{subsec:reviewRobustness}

We briefly review concepts and results needed in the next section.

\begin{definition}
	A pair of reactions in a power law system is called a \textbf{Shinar-Feinberg pair} (SF-pair) in species $X$ if their kinetic order vectors differ only in $X$. A network that contains an SF-pair is said to be of \textbf{Shinar-Feinberg type} (of SF-type).
\end{definition}

The pioneering work of G. Shinar and M. Feinberg was extended by Fortun et al \cite{FMRL2018} to the following theorem:

\begin{theorem}
	Let $\mathscr{N} = (\mathscr{S}, \mathscr{C}, \mathscr{R})$ be a deficiency one CRN and suppose that $(\mathscr{N}, K)$ is a PL-RDK system which admits a positive equilibrium. If a pair of reactions form an SF-pair in species $X$, then the system has ACR in $X$.
\end{theorem}

Fortun and Mendoza \cite{FOME2021} showed that the presence of SF-pairs in a linkage class implied concentration robustness in weakly reversible deficiency zero PL-RDK systems.

The next proposition enables us to determine if a power law system has ACR in a species. Unlike other similar results on ACR, this proposition does not have any deficiency restriction imposed on the underlying network. This allows us to deal with higher deficiency systems including higher deficiency mass action systems. Thus, we consider it as a framework for constructing systems with ACR in a species using ``building blocks'' with low deficiency ($\delta = 0$ or $\delta = 1$).

\begin{proposition}
	Let $(\mathscr{N}, K)$ be a power law system with a positive equilibrium and an independent decomposition $\mathscr{N} = \mathscr{N}_1 \cup \ldots \cup \mathscr{N}_k$. If there is a subnetwork $(\mathscr{N}_i, K_i)$ of deficiency $\delta_i$ with SF-pair in species $X$ such that
	\begin{enumerate}[(i)]
		\item $\delta_i = 0$ and $(\mathscr{N}_i, K_i)$ is a weakly reversible PL-RDK system with the SF-pair in a linkage class; or
		\item $\delta_i = 1$, $(\mathscr{N}_i, K_i)$ is a PL-RDK system, and the SF-pair's reactant complexes are nonterminal;
	\end{enumerate}
	then $(\mathscr{N}, K)$ has ACR in $X$.
\end{proposition}

The analogous ``building blocks'' framework for BCR is the following:

\begin{proposition}
	Let $(\mathscr{N}, K)$ be a power law system with a complex balanced equilibrium and an incidence independent decomposition $\mathscr{N} = \mathscr{N}_1 \cup \ldots \cup \mathscr{N}_k$. If there is a subnetwork $(\mathscr{N}_i, K_i)$ of deficiency $\delta_i$ with SF-pair in species $X$ such that
	\begin{enumerate}[(i)]
		\item $\delta_i = 0$ and $(\mathscr{N}_i, K_i)$ is a weakly reversible PL-RDK system with the SF-pair in a linkage class; or
		\item $\delta_i = 1$, $(\mathscr{N}_i, K_i)$ is a PL-RDK system, and the SF-pair's reactant complexes are nonterminal;
	\end{enumerate}
	then $(\mathscr{N}, K)$ has BCR in $X$.
\end{proposition}

Computational procedures for these frameworks were presented in \cite{FOMF2021}.

\subsection{Concentration Robustness in a Class of LP PL-RDK Systems of SF-Type}
\label{subsec:plrdk}

The following proposition is the basis for broadening the frameworks:

\begin{proposition}
	\label{prop:5.4}
	Let $(\mathscr{N}, K)$ be a cycle terminal PL-RDK system of LP type with flux subspace $\tilde{S}$. If the system has an SF-pair in species $X$ in a linkage class, then it has ACR in $X$.
\end{proposition}

\textit{Proof}. Without loss of generality, we consider the PLP case. If $\tilde{y}$ and $\tilde{y}'$ constitute the kinetic complexes of the SF-pair, then $\tilde{y} - \tilde{y}'$ is an element of $\tilde{S}$. By definition, $\tilde{y} - \tilde{y}'$ has zero coordinates except in $X$. Hence, the scalar product $\langle v, \tilde{y} - \tilde{y}' \rangle = 0$ if and only if $v_X = 0$. Hence, if $v \in \tilde{S}^\perp$, then $v_X = 0$. By the species hyperplane criterion, the system has ACR in $X$. \hfill $\Box$

\begin{remark}
	Proposition \ref{prop:5.4} generalizes an earlier result of Jose et al \cite{JMT2021} on weakly reversible PL-TIK systems with consecutive SF-pairs.
\end{remark}

The broadened framework for ACR is the following:

\begin{proposition}
	\label{prop:LPofSFType}
	Let $(\mathscr{N}, K)$ be a power law system with a positive equilibrium and an independent decomposition $\mathscr{N} = \mathscr{N}_1 \cup \ldots \cup \mathscr{N}_k$. If there is a subnetwork $(\mathscr{N}_i, K_i)$ of deficiency $\delta_i$ with SF-pair in species $X$ such that
	\begin{enumerate}[(i)]
		\item $\delta_i = 0$ and $(\mathscr{N}_i, K_i)$ is a weakly reversible PL-RDK system with the SF-pair in a linkage class; or
		\item $\delta_i = 1$, $(\mathscr{N}_i, K_i)$ is a PL-RDK system, and the SF-pair's reactant complexes are nonterminal; or
		\item $(\mathscr{N}_i, K_i)$ is cycle terminal, a PL-RDK system with the SF-pair in a linkage class, and a PLP system with flux subspace $\tilde{S}_i$;
	\end{enumerate}
	then $(\mathscr{N}, K)$ has ACR in $X$.
\end{proposition}


\section{Concentration Robustness in Poly-PL Systems of LP Type}
\label{sec:concenRobust}

In this section, we identify subsets of poly-PL systems which are of LP type to which we can apply our new results. We first review the known results about LP poly-PL systems and then derive sufficient conditions for the subsets of PL-equilibrated and PL-complex balanced poly-PL systems to be PLP and CLP, respectively. We verify that weakly reversible poly-power law systems with $\hat{T}$-rank maximal kinetics (PY-TIK) are PL-complex balanced and are CLP systems. Towards the end of the section, we establish concentration robustness with an example of an evolutionary game with replicator dynamics.

\subsection{A Review of Previous Results on Poly-PL Systems of LP Type}
\label{subsec:reviewPYK}

The set of poly-PL kinetics was introduced in 2019 by Talabis et al \cite{TMMNJ2020} where they arose in reaction network representations of evolutionary games and was further studied by Magpantay et al in \cite{MHDMN2020}.

\begin{definition}
	\label{def:PYK}
	A kinetics $K: \mathbb{R}^\mathscr{S}_> \rightarrow \mathbb{R}^\mathscr{R}$ is a \textbf{poly-PL kinetics} (PYK) if
	\begin{equation*}
		K_i (x) = k_i \sum_{j=1}^{h_i} a_{ij} x^{f_{ij}} \text{ for all } i = 1, \ldots, r
	\end{equation*}
	written in lexicographic order with $k_i > 0$, $a_{ij} \geq 0$, $f_{ij} \in \mathbb{R}$, and $h_i$ the number of terms in reaction $i$. If $h = \max h_i$, we normalize the length of each kinetics to $h$ by replacing the last term with $h - h_i + 1$ terms with $\displaystyle \frac{1}{h - h_i + 1} x^{f_{i,h_i}}$. We call this the \textbf{canonical PL-representation} of a poly-PL kinetics. We refer to a CRN with PYK as a \textbf{poly-PL system}.
\end{definition}

We set $K_j (x) := k_{ij} x^{f_{ij}}$ with $k_{ij} = k_i a_{ij}$ where $i = 1, \ldots, r$ for each $j = 1, \ldots, h$. $(\mathscr{N}, K_j)$ is a power law system with kinetic order matrix $F_j$.

From Definition \ref{def:PYK}, a PYK can be represented (in a canonical manner) as the sum of PLK. This allows the extension of various results on power law systems to interesting subsets of poly-PL systems. In \cite{TMMNJ2020}, the set of PY-TIK systems is shown to coincide with the set of sums of PL-TIK systems and that any weakly reversible PY-TIK system is complex balanced for all rate constants. Furthermore, it is shown that for a complex balanced PY-TIK system, $Z_+ (\mathscr{N}, K) = \{ x \in \mathbb{R}^\mathscr{S}_> \mid \log x - \log x^* \in \tilde{S}_j^\perp \}$ where $\tilde{S}_j$ is the kinetic reactant flux subspace of any of the $(\mathscr{N}, K_j)$.

\subsection{PL-Complex Balanced and PL-Equilibrated Poly-PL Kinetics}
\label{subsec:PYK}

We recall the definitions of two subsets of complex balanced and positively equilibrated poly-PL systems.

\begin{definition}
	A poly-PL system $(\mathscr{N}, K)$ is \textbf{PL-complex balanced} if the set of complex balanced equilibria $\displaystyle \varnothing \neq Z_+ (\mathscr{N}, K) = \bigcap_{j=1}^h Z_+ (\mathscr{N}, K_j)$. Analogously, $(\mathscr{N}, K)$ is \linebreak \textbf{PL-equilibrated} if the set of positive equilibria $\displaystyle \varnothing \neq E_+ (\mathscr{N}, K) = \bigcap_{j=1}^h E_+ (\mathscr{N}, K_j)$.
\end{definition}

\begin{example}
	As mentioned above, it has been shown that any weakly reversible \linebreak PY-TIK system is unconditionally complex balanced, i.e., it is complex balanced for any set of rate constants. Furthermore, for each such system, $Z_+ (\mathscr{N}, K) = Z_+ (\mathscr{N}, K_j)$ for each $j = 1, \ldots, h$. This implies that $\displaystyle Z_+ (\mathscr{N}, K) = \bigcap_{j=1}^h Z_+ (\mathscr{N}, K_j)$, i.e., any PY-TIK system is PL-complex balanced. It is also clearly a CLP system.
\end{example}

The following proposition generalizes the above example:

\begin{proposition}
	Let $(\mathscr{N}, K)$ be a PL-complex balanced poly-PL system with $(\mathscr{N}, K_j)$ of CLP type with flux subspace $P_{Z,j}$ for $j = 1, \ldots, h$. Then $(\mathscr{N}, K)$ is a CLP system with flux subspace $\displaystyle P_Z = \sum_{j=1}^h P_{Z,j}$.
\end{proposition}

\textit{Proof}. Clearly, $\displaystyle \bigcap_{j=1}^h Z_+ (\mathscr{N}, K_j) = \left\{ x \in \mathbb{R}^\mathscr{S}_> \biggm\vert \log x - \log x^* \in \bigcap_{j=1}^h P_{Z,j}^\perp \right\}$. Since we have $\displaystyle \bigcap_{j=1}^h P_{Z,j}^\perp = \left( \sum_{j=1}^h P_{Z,j} \right)^\perp$ and, by assumption $\displaystyle \varnothing \neq Z_+ (\mathscr{N}, K) = \bigcap_{j=1}^h Z_+ (\mathscr{N}, K_j)$, we obtain the claim. \hfill $\Box$

\begin{example}
	After S. M\"{u}ller and G. Regensburger, every complex balanced PL-RDK system is a CLP system. Since a poly-PL system has reactant-determined kinetics \linebreak (PY-RDK) if and only if each  $(\mathscr{N}, K_j)$ for $j = 1, \ldots, h$ is a PL-RDK system, it follows that any PL-complex balanced PY-RDK system is a CLP system.
\end{example}

The following proposition is the analogue for PLP systems:

\begin{proposition}
	Let $(\mathscr{N}, K)$ be a PL-equilibrated poly-PL system with $(\mathscr{N}, K_j)$ of PLP type with flux subspace $P_{E,j}$ for $j = 1, \ldots, h$. Then $(\mathscr{N}, K)$ is a PLP system with flux subspace $\displaystyle P_E = \sum_{j=1}^h P_{E,j}$.
\end{proposition}

\subsection{Concentration Robustness in Evolutionary Games with Replicator Dynamics}
\label{subsec:games}

\subsubsection{Review of Previous Results}
\label{subsubsec:review}

Recall that Talabis et al showed in \cite{TMMNJ2020} that any evolutionary game with replicator dynamics can be represented as a kinetic system $(\mathscr{N}, K)$ where $\mathscr{N}$ consists of $m$ species $X_i$ and $2m$ reactions $\{ X_i \rightarrow 2 X_i, 2 X_i \rightarrow X_i \}$ for $i = 1, \ldots, m$. The kinetics for the $i$th forward reaction is $X_i f_i(x)$ and the kinetics for the $i$th backward reaction is $X_i \phi(x)$ where $f_i(x)$ is the $i$th payoff function and $\displaystyle \phi(x) = \sum_{i=1}^r X_i f_i(x)$ is the average payoff. Note that the representation is only possible if all the payoff functions are nonnegative. Hence, in a poly-PL replicator game, all the coefficients of the payoff functions must be nonnegative.

The main result in \cite{TMMNJ2020} about poly-PL replicator games is the following:

\begin{theorem}
\label{thm:replicator}
	($m$ variables, $h'$ terms) If the payoff functions of an $m$-variable replicator system with $h'$ terms are of the form
	\begin{equation*}
		f_p(x) = \sum_{i=1}^{h'} \left( a_{pi} \prod_{j=1}^m X_j^{g_{ij}^p} \right) \text{ where } 1 \leq p \leq m,
	\end{equation*}
	if for each $j$ the sets $G^{ij} = \{ g_{ij}^p \mid 1 \leq p \leq m, p \neq j \}$ for $1 \leq j \leq m$ where $1 \leq i \leq h'$ are singleton $\{ g^{ij} \}$ such that $g^{ij} \neq g_{ij}^j$ for $1 \leq j \leq m$ where $1 \leq i \leq h'$, then the replicator system has a positive equilibrium (necessarily complex balanced).
\end{theorem}

The sufficient condition ensures that the associated augmented $T$-matrices have maximal column rank for each $j$ and hence the system is a weakly reversible PY-TIK system which is then unconditionally complex balanced, i.e., has a complex balanced equilibrium for any set of rate constants (Theorem 3 of \cite{TMMNJ2020}).

\subsubsection{The Uniqueness of the Positive Equilibrium of a Poly-PL Replicator Game}
\label{subsubsec:uniqueness}

The result in the paper that implies the uniqueness of the equilibrium is Theorem 4 of \cite{TMMNJ2020}:

\begin{theorem}
\label{thm:uniqueness}
	Let $(\mathscr{N}, K)$ be a weakly reversible poly-PL system with poly $T$-matrices $T_1, \ldots, T_h$ and kinetic reactant deficiency $\hat{\delta} = 0$. Consider an arbitrary poly $T$-matrix $T_k$.
	\begin{enumerate}[(i)]
		\item If $Z_+ (\mathscr{N}, K) \neq \varnothing$ and $x^* \in Z_+ (\mathscr{N}, K)$, then
					\begin{equation*}
						Z_+ (\mathscr{N}, K) = \{ x \in \mathbb{R}^m_> \mid \log x - \log x^* \in \tilde{S}_k^\perp \}
					\end{equation*}
					where $\tilde{S}_k$ is the kinetic reactant flux subspace of $(\mathscr{N}, K_k)$.
		\item If $Z_+ (\mathscr{N}, K) \neq \varnothing$, then $\vert Z_+ (\mathscr{N}, K) \cap Q_k \vert = 1$ for each positive kinetic reactant flux class $Q_k$.
	\end{enumerate}
\end{theorem}

\begin{remark}
	The index $k$ in Theorem \ref{thm:uniqueness} is the index $i$ in Theorem \ref{thm:replicator}.
\end{remark}

We formulate our result as a corollary:

\begin{corollary}
	\label{cor:6.9}
	(Corollary to Theorem \ref{thm:replicator}) Any poly-PL replicator game satisfying the sufficient condition of Theorem \ref{thm:replicator} has a unique positive equilibrium.
\end{corollary}

\textit{Proof}. It follows directly from the definition of the block matrix as direct sum of the poly $T$-matrices that a PY-RDK system $(\mathscr{N}, K)$ is a PY-TIK system if and only if each of its summands $(\mathscr{N}, K_j)$ for $j = 1, \ldots, h$ is a PL-TIK system. Since after \cite{BOROS2012}, any weakly reversible PL-TIK system is unconditionally complex balanced, it follows from the result of S. M\"{u}ller and G. Regensburger that such a PL-TIK system has zero kinetic deficiency, i.e., $\tilde{\delta} = 0$. This is equivalent to $\tilde{s}_j = n_j - \ell$ where $\tilde{s}_i$ is the rank of $(\mathscr{N}, K_j)$, $n_j$ is the number of complexes in $(\mathscr{N}, K_j)$, and $\ell$ is the number of linkage classes of $\mathscr{N}$. For the replicator system, $n - \ell = 2m - m = m$ where $n$ and $m$ represent the number of complexes and species, respectively, of the network. Thus, the kinetic order subspace $\tilde{S}_j = \mathbb{R}^m$. Hence $Z_+ (\mathscr{N}, K_j) = \{ x \in \mathbb{R}^m_> \mid \log x - \log x^* \in \tilde{S}_j^\perp = 0 \}$, i.e., it consists of a single element. According to $(i)$ of Theorem \ref{thm:uniqueness}, $Z_+ (\mathscr{N}, K)$ coincides with this set. \hfill $\Box$

\begin{remark}
	\textcolor{white}{}
	\begin{enumerate}[(i)]
		\item The conclusion of Corollary \ref{cor:6.9} also follows from $(ii)$ of Theorem \ref{thm:uniqueness}, which is a general consequence for a CLP system.
		\item Corollary \ref{cor:6.9} can be generalized to the following proposition: Let $(\mathscr{N}, K)$ be a weakly reversible PL-RDK system with zero kinetic deficiency and $n - \ell = m$ where $m$, $n$, and $\ell$ represent the number of species, number of complexes, and number of linkage classes, respectively. Then $(\mathscr{N}, K)$ has a unique complex balanced equilibrium for each set of rate constants.
		\item In \cite{BOROS2013}, it is shown that any PL-complex balanced PY-RDK system is a CLP system with parameter subspace $\displaystyle \left( \sum_{j=1}^h \tilde{S}_j \right)^\perp$ where the $\tilde{S}_j$ are the kinetic order subspaces (after S. M\"{u}ller and G. Regensburger) of the PL-RDK summands of the PY-RDK system. In a PY-TIK system, the intersection is over the same parameter subspace. In other words, the parameter subspace of a PY-TIK system is smaller than that of other PL-complex balanced PY-RDK systems, resulting in a larger equilibria parameter subspace.
	\end{enumerate}
\end{remark}

\section{Summary and Conclusions}
\label{sec:summary}

A CKS $(\mathscr{N}, K)$ is a PLP system if its set of positive equilibria is an LP set, i.e., $E_+ (\mathscr{N}, K) = E(P_E, x^*)$ where $P_E$ is the flux subspace and $x^*$ is a given positive equilibrium. Analogously, the CKS is a CLP system if its set of complex balanced equilibria $Z_+ (\mathscr{N}, K) = E(P_Z, x^*)$ where $P_Z$ is the flux subspace and $x^*$ is a given complex balanced equilibrium. Various LP systems beyond mass action systems have already been studied. In this paper, we studied the properties of concentration robustness of species in LP systems. 

We presented the species hyperplane criterion, a necessary and sufficient condition for ACR and BCR. Our approach was based on a necessary and sufficient condition for ACR in PLP systems and its analogue for BCR in CLP systems in terms of their LP subspaces.  In Theorem \ref{thm:hyperplanecriterion}, we presented a necessary condition yielding an upper bound for the number of species in which the system admits concentration robustness. Through Propositions \ref{prop:4.1} and \ref{prop:4.4}, it can be determined in which species concentration robustness in an LP system holds and which LP systems can be used as control components for concentration robustness in decomposition subnetworks. Further application of the species hyperplane criterion for concentration robustness in LP systems of SF-type (presented in Proposition \ref{prop:LPofSFType}) led us to broaden the ``building blocks'' framework for constructing larger systems with concentration robustness studied in \cite{FOMF2021, FOME2021}. It allows easier identification of concentration robustness in larger networks with arbitrary deficiency. The results were applied to species concentration robustness in LP systems with poly-PL kinetics from evolutionary game theory (see Section \ref{sec:concenRobust}).





\appendix

\section{List of Abbreviations}

\begin{table}[ht]
    \begin{center}
    \caption{List of abbreviations}\label{tab1}
    \begin{tabular}{@{}ll@{}}
        \hline
		Abbreviation & Meaning \\
		\hline
		ACR & Absolute Concentration Robustness \\
		BCR & Balanced Concentration Robustness \\
		CKS & Chemical Kinetic System \\
		CLP & Complex Balanced Equilibria Log-Parametrized \\
		CRN & Chemical Reaction Network \\
		GMAS & Generalized Mass Action System \\
		LP & Log-Parametrized \\
		MAK & Mass Action Kinetics \\
		ODE & Ordinary Differential Equation \\
		PL-NDK & Power Law System with Non-Reactant-Determined Kinetics \\
		PL-RDK & Power Law System with Reactant-Determined Kinetics \\
		PL-TIK & Power Law System with $\hat{T}$-Rank Maximal Kinetics \\
		PLK & Power Law Kinetics \\
		PLP & Positive Equilibria Log-Parametrized \\
		PY-RDK & Poly-Power Law System with Reactant-Determined Kinetics \\
		PY-TIK & Poly-Power Law System with $\hat{T}$-Rank Maximal Kinetics \\
		PYK & Poly-Power Law Kinetics \\
		SF & Shinar-Feinberg \\
		SFRF & Species Formation Rate Function \\
		\hline
    \end{tabular}
    \end{center}
\end{table}

\end{document}